\documentclass[12pt,oneside]{amsart}
\usepackage{amsmath}
\usepackage{amsthm}
\usepackage{amsfonts}
\usepackage{a4wide}

\usepackage{amssymb}
\newtheorem{thm}{Theorem}
\newtheorem{lemma}[thm]{Lemma}
\newtheorem{prop}[thm]{Proposition}
\newtheorem{corollary}[thm]{Corollary}

\theoremstyle{definition}

\theoremstyle{remark}
\newtheorem{rem}[thm]{Remark}
\newcommand{\DD}{{\mathbb D}}
\newcommand{\OO}{{\mathcal O}}

\newcommand{\EE}{{\mathbb E}}
\newcommand{\GG}{{\mathbb G}}

\newcommand{\RR}{{\mathbb R}}
\newcommand{\CC}{{\mathbb C}}

\DeclareMathOperator{\Aut}{Aut} \DeclareMathOperator{\id}{id}
 
 \DeclareMathOperator{\re}{Re}

\renewcommand{\phi}{\varphi}

\begin{document}



\title{Schwarz lemma for the tetrablock}

\author{Armen Edigarian}

\address{Institute of Mathematics, Jagiellonian University,
Reymonta 4, 30-059 Krak\'ow, Poland}
\email{{Armen.Edigarian}@im.uj.edu.pl}
\thanks{The paper was partially supported by the Research Grant of the Polish Ministry of Science and Higher Education
No. 1 PO3A 005 28.}

\author{W\l odzimierz Zwonek}

\address{Institute of Mathematics, Jagiellonian University,
Reymonta 4, 30-059 Krak\'ow, Poland}
\email{{Wlodzimierz.Zwonek}@im.uj.edu.pl}





\begin{abstract} We describe all complex geodesics in the tetrablock passing through the origin thus obtaining the form of all extremals
in the Schwarz Lemma for the tetrablock. Some other extremals for
the Lempert function and geodesics are also given. The paper may
be seen as a continuation of the results from \cite{AWY}. The
proofs rely on a necessary form of complex geodesics in general
domains which is also proven in the paper.

\end{abstract}

\maketitle

\section{Introduction and main results}

Let us recall two holomorphically invariant families of functions.
For a domain $D\subset\CC^n$, $w,z\in D$ define the {\it Lempert
function $\tilde k_D$} as follows
\begin{equation*}
\tilde k_D(w,z):=\inf\{p(\lambda_1,\lambda_2):\text{ there is
$f\in\OO(\DD,D)$ such that $f(\lambda_1)=w$, $f(\lambda_2)=z$}\}
\end{equation*}
and the \emph{Carath\'eodory pseudodistance} $c_D$
\begin{equation*}
c_D(w,z):=\sup\{p(F(w),F(z)):F\in\OO(D,\DD)\},
\end{equation*}
where $\OO(\Omega_1,\Omega_2)$ denotes the set of all holomorphic
mappings $\Omega_1\to \Omega_2$ and $p$ denotes the Poincar\'e distance on the unit disc $\DD$.
Basic properties and results on
both functions as well as definitions of taut and hyperconvex
domains that we use in our paper one may find in \cite{J-P}.

It is well-known that $c_D\leq\tilde k_D$ on any domain $D$.


In 1981 L. Lempert proved the following important result (see \cite{Lem 1981}, \cite{J-P}, \cite{Edi1}).
\begin{thm}\label{thm:0}
Let $D$ be a domain in $\CC^n$ and let $\{D_j\}_{j=1}^\infty$ be a
sequence of domains exhausting $D$, i.e., $D_1\subset
D_2\subset\dots\subset D$ and $\cup_{j=1}^\infty D_j=D$. Assume
that any $D_j$ is biholomorphic to a convex domain. Then
$c_D=\tilde k_D$.
\end{thm}
For more then twenty years in was an open question whether any
bounded pseudoconvex domain $D\subset\CC^n$ such that $c_D=\tilde
k_D$ can be exhausted by domains biholomorphically equivalent to a
convex domain.

In 2001 J. Agler and N.J. Young introduced the symmetrized bidisc
$\GG_2$ (see \cite{A-Y1}), a domain which turns up in control
engineering and produces problems of a function-theoretic
character. In particular, the equality $c_{\GG_2}(0,\cdot)=\tilde
k_{\GG_2}(0,\cdot)$ on $\GG_2$ has been shown in \cite{A-Y1}. This
result was extended later to the identity $c_{\GG_2}\equiv\tilde
k_{\GG_2}$ (see \cite{Cos 2004}, \cite{Agl-You 2004}). The
symmetrized bidisc and its higher dimensional analogue, the
symmetrized polydisc $\GG_n$, has been recently extensively
studied (see e.g. \cite{A-Y1}, \cite{Agl-You 2004}, \cite{Cos
2004}, \cite{Pfl-Zwo 2005}, \cite{Edi-Zwo 2005}, \cite{Agl-You
2006} and others). It was shown that $\GG_2$ cannot be exhausted
by domains biholomorphic to convex ones (see \cite{Cos 2004},
\cite{Edi2}). So symmetrized bidisc delivers a counterexample to
the above problem.

In two papers \cite{AWY} and \cite{You 2007} and in the PhD Thesis
\cite{Abo 2007} the Authors initiated the study of another domain
which also naturally appears in control engineering and produces
problems of a function-theoretic character. The {\it tetrablock}
is a domain in $\CC^3$, denoted by $\EE$, which is the image of
the Cartan domain $R_I(2,2)$ upon the mapping
$\pi(A):=(a_{11},a_{22},\det A)$, $A=(a_{jk})_{j,k=1,2}\in\CC^{2\times 2}$. In the paper \cite{AWY} several
equivalent definitions of the domain $\EE$ are given. Recall one
of them
\begin{equation}\label{eq:tetrablock}
\EE=\{(z_1,z_2,z_3)\in\CC^3: |z_1-\bar z_2z_3|+|z_2-\bar z_1z_3|+|z_3|^2<1\}.
\end{equation}

It is proven in \cite{AWY} that the equality between
the Carath\'eodory distance and the Lempert function of $\EE$ with
fixed at the origin one of the arguments,
$$
c_{\EE}(0,\cdot)=\tilde k_{\EE}(0,\cdot)
$$
holds on $\EE$, which suggests that the equality between both
functions could hold on $\EE\times\EE$. In our paper we deal with
this domain. We find all the solutions of the extremal problem in
the Schwarz Lemma for the tetrablock (in other words all functions
$f\in\OO(\DD,\EE)$ such that $f(0)=0$ and $\tilde
k_{\EE}(f(0),f(\lambda))=p(0,\lambda)$, $\lambda\in\DD$).
Moreover, we find some other $\tilde k_{\EE}$-extremals and
complex geodesics and we make several observations which may lead
to the solution of the problem whether the Lempert Theorem is
valid on $\EE$.


The equality of the Lempert function and the Caratheodory
pseudodistance on a given domain is closely related to the notion
of a complex geodesic. A holomorphic mapping $f:\DD\to D$ is
called a \emph{(complex) geodesic} if there exists a mapping
$F:D\to\DD$ such that $F\circ f$ is an automorphism of $\DD$
(without loss of generality, one can assume that $F\circ
f=\id_{\DD}$). Note that if $f:\DD\to D$ is a complex geodesic
then $c_D=\tilde k_D$ on $f(\DD)\times f(\DD)$.

It is proven in \cite{AWY} that on the tetrablock complex geodesics passing through
the origin and another, arbitrary, point of $\EE$ exist. Moreover, for any $z\in \EE$ a geodesic
passing through $0$ and $z$ is given. In our paper we give a
description of all such geodesics. This result may be seen as a
Schwarz-type Lemma for the tetrablock.

\begin{thm}\label{thm:2} Let $\phi:\DD\to\bar\DD$ be a holomorphic
mapping such that $\phi(0)=-C$, where $C\in[0,1]$. Then for any
$\omega_1,\omega_2\in\partial\DD$ the mapping
\begin{equation}\label{eq:2}
f(\lambda)=\Big(\omega_1\frac{\phi(\lambda)+C}{1+C},\omega_2\lambda\frac{1+C\phi(\lambda)}
{1+C},\omega_1\omega_2\lambda\phi(\lambda)\Big)
\end{equation}
is a complex geodesic in $\EE$. Moreover, any complex geodesic $f:\DD\to\EE$ such that $f(0)=0$
is (up to a permutation of two first variables)
of type \eqref{eq:2}.
\end{thm}



Since the set $S:=\{z\in\EE: z_1z_2=z_3\}$ is the orbit of $0$ of the group $\Aut \EE$
(see \cite{AWY} and \cite{You 2007}), the above theorem gives the full description of complex geodesics intersecting
$S$.

In the proof of Theorem~\ref{thm:2} the main role is played by a
result giving a necessary condition on complex geodesics in
domains with a lot of automorphisms. Similar characterizations one
can find e.g. in \cite{Pol}, \cite{Edi}, see also \cite{J-P}.
\begin{thm}\label{thm:10} Let $D\subset\CC^n$ be a domain and let $f:\DD\to D$, $F:D\to\DD$ be holomorphic mappings
such that $F\circ f=\id_{\DD}$. Suppose that
$\Phi:(-\tau,\tau)\times D\to D$, $\tau>0$, is a $C^2$ mapping
with
\begin{enumerate}
\item $\Phi_0=\id$, where $\Phi_t=\Phi(t,\cdot)$;
\item $\Phi_t:D\to D$ is a holomorphic mapping.
\end{enumerate}
Then there exists a $C\in\RR$ such that
\begin{equation*}
\psi(\lambda)=-\bar\psi(0)\lambda^2+iC\lambda+\psi(0)\quad\text{ for any }\lambda\in\DD,
\end{equation*}
where
\begin{equation*}
\psi(\lambda)=\sum_{j=1}^n \frac{\partial F}{\partial z_j}(f(\lambda))\gamma_j(f(\lambda))
\end{equation*}
and $\gamma=(\gamma_1,\dots,\gamma_n):D\to\CC^n$ is a vector field on $D$ generated by a family $\{\Phi_t\}_{|t|<\tau}$,
i.e., $\gamma(z)=\lim_{t\to0}\frac{\Phi_t(z)-z}{t}$.
\end{thm}

In the paper we also deliver some other $\tilde k_{\EE}$-extremals
(which, assuming the equality $\tilde k_{\EE}=c_{\EE}$ would be
complex geodesics) which are essentially different from the ones
given in Theorem~\ref{thm:2} (see Corollary~\ref{cor:3}); we make
also some observations which may lead to the better understanding
of the geometry of $\EE$.

\section{Circular domains and complex geodesics}

We start this section with the proof of Theorem~\ref{thm:10}.

\begin{proof} Fix $\lambda\in\DD\setminus\{0\}$. Consider functions $a(t)=F(\Phi_t(f(\lambda)))$
and $b(t)=F(\Phi_t(f(0)))$. Note that $a(0)=\lambda$, $a'(0)=\psi(\lambda)$, $b(0)=0$ and $b'(0)=\psi(0)$.

>From the Schwarz-Pick lemma (see e.g. \cite{J-P}) for any $t\in(-\tau,\tau)$ we have
\begin{equation*}
\Big|\frac{a(t)-b(t)}{1-\overline b(t) a(t)}\Big|\le|\lambda|
=\Big|\frac{a(0)-b(0)}{1-\overline b(0) a(0)}\Big|.
\end{equation*}
Therefore, the function
$\rho(t):=\Big|\frac{a(t)-b(t)}{1-\overline b(t) a(t)}\Big|^2$
attains its maximum at $t=0$ and, therefore, $\rho'(0)=0$. We have
\begin{equation*}
\rho'(0)=2\re\Big(a'(0)\overline a(0)-\overline a(0)b'(0)(1-|a(0)|)^2\Big).
\end{equation*}
From this we get
$\re\Big(\frac{\psi(\lambda)-\psi(0)}{\lambda}+\lambda\overline\psi(0)\Big)=0$
for any $\lambda\in\DD\setminus\{0\}$. Since $\psi$ is a
holomorphic function on $\DD$ there exists a constant $C\in\RR$
such that
$\frac{\psi(\lambda)-\psi(0)}{\lambda}+\lambda\overline\psi(0)=iC$.
\end{proof}

Before we apply Theorem~\ref{thm:10} in the proof of
Theorem~\ref{thm:2} we make some remarks how this theorem may be
applied in the study of complex geodesics in domains with many
symmetries.

The families $\{\Phi_t\}$ we are going to consider are
one-parameter groups of transformations on circular type domains.
Recall that $D\subset\CC^n$ is {\it
$(\alpha_1,\dots,\alpha_n)$-circular} if $(e^{i\alpha_1
t}z_1,\dots,e^{i\alpha_n t}z_n)\in D$ for any $(z_1,\dots,z_n)\in
D$ and any $t\in\RR$. Here $\alpha_1,\dots,\alpha_n\in\RR$. For
any such a domain we consider a mapping $\Phi:\RR\times D\to D$
defined as
$\Phi(t,z)=(e^{it\alpha_1}z_1,\dots,e^{it\alpha_n}z_n)$. Note that
in this case $\gamma(z)=i(\alpha_1 z_1,\dots, \alpha_n z_n)$. So,
from Theorem~\ref{thm:10} we get
\begin{corollary}\label{cor:2}
Let $D\subset\CC^n$ be an $(\alpha_1,\dots,\alpha_n)$-circular domain and let $f:\DD\to D$,
$F:D\to\DD$ be holomorphic mappings  such that $F\circ f=\id$. Then there exist constants $C\in\RR$, $a\in\CC$ such that
\begin{equation*}
\sum_{j=1}^n \alpha_j\frac{\partial F}{\partial z_j}(f(\lambda))f_j(\lambda)=
\overline a\lambda^2+C\lambda+a\quad\text{ for any }\lambda\in\DD.
\end{equation*}
\end{corollary}
Now we are giving some examples of applications of our result. Let
us start with Reinhardt domains, which are seemingly well
understood (see \cite{Zwo}). Even in this case we get new results.

\subsection{Reinhardt domains} Recal that a domain $D\subset\CC^n$ is called {\it Reinhardt} if it is
$(\alpha_1,\dots,\alpha_n)$-circular for any $\alpha_1,\dots,\alpha_n\in\RR$.

\begin{corollary} Let $D\subset\CC^n$ be a Reinhardt domain and let $f:\DD\to
D$, $F:D\to\DD$ be holomorphic mappings such that $F\circ
f=\id_{\DD}$. Then for any $j=1,\dots,n$ there exist constants
$a_j\in\CC$ and $C_j\in\RR$ such that
\begin{equation*}
\frac{\partial F}{\partial z_j}(f(\lambda))f_j(\lambda)=\overline a_j\lambda^2+C_j\lambda+a_j\quad
\text{ for any }\lambda\in\DD.
\end{equation*}
In particular, for any complex geodesic $f=(f_1,\dots,f_n):\DD\to
D$ and for any $j=1,\dots,n$ the following situation holds: $f_j$
has at most one zero in $\DD$, $f_j\equiv0$ or $\frac{\partial
F}{\partial z_j}(f(\cdot))\equiv 0$.
\end{corollary}
Note that the example of the bidisc and geodesics of the form
$f:\DD\owns\lambda\to(\lambda,b(\lambda))\in\DD^2$ with
$b\in\OO(\DD,\DD)$, $b(0)=0$ shows that there are geodesics such
that $f_j$ has (even infinitely) many zeros but $f_j\not\equiv 0$.



\subsection{The symmetrized bidisc} We call a domain
\begin{equation*}
\GG_2=\{(\lambda+\mu,\lambda\mu): \lambda,\mu\in\DD\}
\end{equation*}
the {\it symmetrized bidisc} (see \cite{A-Y1}). Note that the
symmetrized bidisc is $(1,2)$-circular. The uniqueness result on
geodesics and their  complete description in $\GG_2$ are given in
\cite{Pfl-Zwo 2005} and \cite{Agl-You 2006}. Here, we show how one
can use our results to give a simpler proof of the description of
all geodesics in $\GG_2$ passing through the origin.

Assume that $f=(f_1,f_2):\DD\to\GG_2$ is a complex geodesic such
that $f(0)=0$. It is shown in \cite{A-Y1} that there exists an
$\omega\in\partial\DD$ such that $F\circ f=\id_{\DD}$, where
$F(z_1,z_2)=\frac{2\omega z_2-z_1}{2-\omega z_1}$. Note that
\begin{equation*}
\frac{\partial F}{\partial z_1}(f)f_1=
\frac{-2+2\omega^2 f_2}{(2-\omega f_1)^2}f_1\quad\text{ and }\quad
\frac{\partial F}{\partial z_2}(f)f_2=\frac{2\omega}{2-\omega f_1}f_2.
\end{equation*}
and, therefore,
\begin{equation*}
\frac{-2+2\omega^2 f_2(\lambda)}{(2-\omega f_1(\lambda))^2}f_1+2\frac{2\omega}{2-\omega f_1(\lambda)}f_2(\lambda)=C\lambda
\end{equation*}
for some $C\in\RR$. Since $\frac{2\omega f_2-f_1}{2-\omega f_1}$
is the identity, we get the equality $2\omega
f_2=f_1+\lambda(2-\omega f_1)$ on $\DD$, so the following equality
has to be satisfied on $\DD$
$$
f_1^2(\lambda\omega^2-\omega-C\lambda\omega^2)+f_1(2-6\lambda\omega+4C\lambda\omega)+4\lambda(2-C)=0.
$$
Consequently,
$$
f_1(\lambda)=\frac{-(2-6\lambda\omega+4C\lambda\omega)\pm2(\lambda\omega+1)}{2\omega(\lambda\omega-1-C\lambda\omega)},\;\lambda\in\DD.
$$
Since $f_1(0)=0$ we get that
$f_1(\lambda)=\frac{2(2-C)\lambda}{\omega\lambda(1-C)-1}$, and
then
$f_2(\lambda)=\lambda\frac{\lambda-\bar\omega(1-C)}{1-\lambda\omega(1-C)}$,
$\lambda\in\DD$. Since $f_1(\DD)\subset2\DD$ and
$f_2(\DD)\subset\DD$ applying the Schwarz Lemma we get that
$|2-C|\leq 1$ and $|1-C|\leq 1$; which gives the condition
$C\in[1,2]$. It is simple to see that this condition is not only
necessary but also sufficient (see e. g. \cite{Pfl-Zwo 2005} and
\cite{Agl-You 2006}).

\section{Description of complex geodesics of the tetrablock passing through the origin}
Before we present the proof of Theorem~\ref{thm:2} let us present
some basic properties of the tetrablock. The properties that we
present below are either obvious or come from \cite{AWY},
\cite{Abo 2007} or \cite{You 2007}.

It is obvious that the permutation of two first variables, i.e.
the mapping $\sigma$ given by the formula
$\sigma(z_1,z_2,z_3)=(z_2,z_1,z_3)$, is an automorphism of $\EE$.

For any $\omega\in\partial\EE$ the mapping $F_{\omega}$ given by
the formula $F_{\omega}(z)=(\omega z_1,z_2,\omega z_3)$, is also
an automorphism of $\EE$. Hence, the tetrablock is $(1,0,1)$ and
$(0,1,1)$-circular.

Recall that $z\in\EE$ iff $\left|\Psi_{\eta}(x)\right|<1$ for any
$\eta\in\bar\DD$, where
\begin{equation*}
\Psi_{\eta}(z_1,z_2,z_3)=\frac{\eta z_3-z_2}{\eta z_1-1}.
\end{equation*}

It is quite simple to see that there is a continuous function
$\rho:\EE\mapsto[0,\infty)$ such that $\log \rho$ is
plurisubharmonic and $\rho(\lambda z_1,\lambda z_2,\lambda^2
z_3)=|\lambda|\rho(z)$ for any $z=(z_1,z_2,z_3)\in\EE$,
$\lambda\in\CC$. This shows that $\EE$ is a hyperconvex domain; in
particular, $\EE$ is a taut domain, too.

Note also that for any $z\in\EE$ the inequalities $|z_j|<1$,
$j=1,2,3$, hold.

\begin{proof}[Proof of Theorem~\ref{thm:2}]
Assume that $\phi,\psi:\DD\to\DD$ are holomorphic mappings and
that $C\in[0,1)$. Then for any $\omega_1,\omega_2\in\partial\DD$
the mapping
\begin{equation}\label{eq:3}
f(\lambda)=\Big(\omega_1\frac{\phi(\lambda)+C}{1+C},\omega_2\psi(\lambda)\frac{1+C\phi(\lambda)}
{1+C},\omega_1\omega_2\phi(\lambda)\psi(\lambda)\Big)
\end{equation}
satisfies the inclusion $f(\DD)\subset\EE$.

Indeed, we have
\begin{align*}
|f_1-\overline f_2 f_3|&=\Big|\frac{\phi+C}{1+C}-\overline\psi\frac{1+C\overline\phi}{1+C}
\phi\psi\Big|=\frac{|\phi(1-|\psi|^2)+C(1-|\phi|^2|\psi|^2)|}{1+C},\\
|f_2-\overline f_1 f_3|&=\Big|\psi\frac{1+\phi
C}{1+C}-\frac{\overline\phi+C}{1+C}
\phi\psi\Big|=\frac{|\psi|(1-|\phi|^2)}{1+C}
\end{align*}
Hence,
\begin{multline*}
|f_1-\overline f_2 f_3|+|f_2-\overline f_1 f_3|+|f_3|^2\le
\frac{|\phi|(1-|\psi|^2)+C(1-|\phi|^2|\psi|^2)}{1+C}+\frac{|\psi|(1-|\phi|^2)}{1+C}\\
+|\phi|^2|\psi|^2<1.
\end{multline*}

Note that $\Psi_{\overline \omega_1}\circ
f(\lambda)=\omega_2\psi(\lambda)$. Hence, if $\psi$ is an
automorphism then $f$ is a complex geodesic. In particular, the
sufficient part of Theorem~\ref{thm:2} is proven (the case $C=1$
and $\phi\equiv -1$, which formally is not covered by the above
reasoning, is simple).

Take now any geodesic $f\in\OO(\DD,\EE)$ such that $f(0)=0$. It
follows from \cite{AWY} there is an $\omega\in\partial\DD$ such
that $\Phi\circ \phi$ is a rotation, where $\Phi=\Psi_{\omega}$ or
$\Phi=\Psi_{\omega}\circ \sigma$.
Since $\sigma$ and $F_{\omega}$ are automorphisms of $\EE$ we
easily arrive at the following statement. In order to show that
all geodesics of $\EE$ with $f(0)=0$ are of the form given in the
theorem it is sufficient to show that all geodesics
$f\in\OO(\DD,\EE)$ with $f(0)=0$ such that $\Psi_{1}\circ f$ are
as in the theorem.

So, take $F=\Psi_1$. We are going to find necessary form of all
right inverses $f:\DD\to\EE$ of $F$ (for a while we assume that
$f(0)$ may take any value in $\EE$). Note that $\frac{\partial
F}{\partial z_1}=\frac{z_2-z_3}{(z_1-1)^2}$ and $\frac{\partial
F}{\partial z_3}=\frac{1}{z_1-1}$. Hence, in view of
Corollary~\ref{cor:2} there exist constants $a\in\CC$ and
$C\in\RR$ such that
\begin{equation}\label{eq:10}
\frac{f_1(\lambda) f_2(\lambda)-f_3(\lambda)}{(f_1(\lambda)-1)^2}=\overline a\lambda^2+C\lambda+a.
\end{equation}
Since $f$ is the right inverse of $F$, we have
$f_3(\lambda)-f_2(\lambda)=\lambda(f_1(\lambda)-1)$,
$\lambda\in\DD$. From this and \eqref{eq:10} we get
\begin{equation*}
\begin{cases}
f_1(\lambda)=&(f_3(\lambda)+\overline a\lambda^2+C\lambda+a)/(\overline a\lambda^2+(C+1)\lambda+a)\\
f_2(\lambda)=&(f_3(\lambda)(\overline a\lambda^2+C\lambda+a)+\lambda^2)/(\overline a\lambda^2+(C+1)\lambda+a).
\end{cases}
\end{equation*}
But $f_1(0)=f_3(0)=0$ so $a=0$ and, therefore,
\begin{equation*}
\begin{cases}
f_1(\lambda)=&(\phi(\lambda)+C)/(C+1)\\
f_2(\lambda)=&\lambda(\phi(\lambda)C+1)/(C+1)\\
f_3(\lambda)=&\lambda \phi(\lambda),
\end{cases}
\end{equation*}
where $\phi:\DD\to\CC$ is a holomorphic mapping. Since
$f_1,f_2,f_3\in\OO(\DD,\DD)$ and $f(0)=0$, we get that
$\phi\in\OO(\DD,\bar\DD)$, $\phi(0)=-C$. The Schwarz Lemma applied
to $f_2$ and $f_3$ gives the inequalities $|C|\leq 1$ and
$|1-C|\leq 1$, so $C\in[0,1]$. Hence, we have finished the proof.
\end{proof}

\begin{rem}
It was shown in \cite{AWY} that no uniqueness of geodesics (even
passing through $0$) is possible in the case of $\EE$. This can
also be easily seen from Theorem~\ref{thm:2}. But we may show also
in a different (simple) way. Namely, consider the following
embedding
\begin{equation*}
\Phi:\DD^2\owns(\lambda,\mu)\mapsto(\lambda,\mu,\lambda\mu)\in\EE.
\end{equation*}
Now the inequalities
\begin{equation*}
d_{\DD^2}((\lambda_1,\mu_1),(\lambda_2,\mu_2))\geq
d_{\EE}((\lambda_1,\mu_1,\lambda_1\mu_1),(\lambda_2,\mu_2,\lambda_2\mu_2))\geq
d_{\DD^2}((\lambda_1,\mu_1),\lambda_2,\mu_2)),
\end{equation*}
where $d=\tilde k$ or $c$, together with the equality $\tilde k_{\DD^2}=c_{\DD^2}$ show that
\begin{equation*}
\tilde
k_{\EE}((\lambda_1,\mu_1,\lambda_1\mu_1),(\lambda_2,\mu_2,\lambda_2\mu_2))=
c_{\EE}((\lambda_1,\mu_1,\lambda_1\mu_1),(\lambda_2,\mu_2,\lambda_2\mu_2))
\end{equation*}
for any $(\lambda_1,\mu_1),(\lambda_2,\mu_2)\in\DD^2$. Moreover,
any mapping of the form
\begin{equation*}
\DD\owns\lambda\mapsto(a(\lambda),b(\lambda),a(\lambda)b(\lambda))\in\EE,
\end{equation*}
where $a,b\in\OO(\DD,\DD)$ and at least one of the functions $a$,
$b$ is an automorphism of $\DD$, is a complex geodesic, which
easily delivers non-uniqueness result for geodesics (passing
through $0$) - to make it visible consider the functions
$\DD\owns\lambda\mapsto(\lambda,b(\lambda),\lambda
b(\lambda))\in\EE$, where $b$ is an arbitrary function (out of
many) from $\OO(\DD,\DD)$ such that $b(0)=0$ and
$b(\frac{1}{2})=\frac{1}{4}$. In fact, it is easy to see that all
the geodesics passing through two points from the set
$S=\{z\in\EE:z_1z_2=z_3\}=\Phi(\DD^2)$ are of the above form. To
see this take two different points of the form
$\Phi(\lambda_1,\mu_1),\Phi(\lambda_2,\mu_2)$ from $\EE$. Without
loss of generality assume that $p(\lambda_1,\lambda_2)\geq
p(\mu_1,\mu_2)$ and the geodesic $f$ is such that
$f(\lambda_1)=\Phi(\lambda_1,\mu_1)$,
$f(\lambda_2)=\Phi(\lambda_2,\mu_2)$. Then
$\Psi_{\omega}\circ\sigma\circ f,f_1\in\OO(\DD,\DD)$ and both
functions map $\lambda_j$ into $\lambda_j$, $j=1,2$, which in view
of the Schwarz-Pick Lemma gives that they are both identities,
which easily shows the desired form of $f$.
\end{rem}

\begin{rem}
Note that the boundary of $\EE$ contains many non-constant holomorphic discs. Fix a constant $C\in[0,1]$ and a
holomorphic function $\phi\in\OO(\DD,\DD)$. Then it follows from calculations made in the proof of Theorem~\ref{thm:2}
that $f:\DD\to\partial\EE$, where
\begin{equation*}
f(\lambda)=\Big(\omega_1\frac{\phi(\lambda)+C}{1+C},\omega_2\frac{1+C\phi(\lambda)}
{1+C},\omega_1\omega_2\phi(\lambda)\Big).
\end{equation*}
\end{rem}

\section{The Lempert function of the tetrablock}
Note that $\EE\cap(\{0\}\times\CC^2)=\{(0,z,w): |z|+|w|<1\}$.
We may calculate the Lempert function for some special points -- note that in the proof
of the result below we make use of $\tilde k_{\EE}$-extremals omitting the special set $S$.

\begin{prop}\label{prop:1} For any $z,w\in\DD$, such that $|z|+|w|<1$ we have
\begin{equation}
\tilde k_{\EE}\big((0,0,w),(0,z,w)\big)=\frac{|z|}{1-|w|}.
\end{equation}
\end{prop}

Recall that a mapping $f\in\OO(\DD,D)$ such that $\tilde
k_D(f(\lambda_1),f(\lambda_2))=p(\lambda_1,\lambda_2)$ for
some $\lambda_1,\lambda_2\in \DD$ is called a {\it $\tilde
k_D$-extremal} for $(f(\lambda_1),f(\lambda_2))$. Recall
that if $D$ is taut then for any $w,z\in D$ there is a $\tilde
k_D$-extremal for $(w,z)$. It is easy to see that if
$f\in\OO(\DD,D)$ is such that $f(\DD)\subset\subset D$ then
$f$ is not a $\tilde k_D$-extremal for any $(w,z)\in D\times D$
with $w\neq z$.


To prove Proposition~\ref{prop:1} we follow the ideas on
transportation of $\tilde k$-extremals (and complex geodesics)
that may be found in \cite{Pfl-Zwo 1998}, \cite{Zwo} (and
references there).

First note that for any $f\in\OO(\DD,\EE)$ with the equality
$f_1(0)=f_3(0)=0$ the function $\tilde f$ given by the
formula
$\tilde f(\lambda):=(\frac{f_1(\lambda)}{\lambda},f_2(\lambda),\frac{f_3(\lambda)}{\lambda})$,
$\lambda\in\DD$ (certainly here and later we mean
$\tilde f_j(0)=f_j^{\prime}(0)$, $j=1,3$), is either from
$\OO(\DD,\partial\EE)$ or from $\OO(\DD,\EE)$. Indeed, to see this apply the maximum principle
for subharmonic functions and the plurisubharmonicity of $\log\rho$.

\begin{lemma}\label{lem:1} Let $f\in\OO(\DD,\EE)$ be a $\tilde
k_{\EE}$-extremal for $(f(0),f(\sigma))$, where
$\sigma\in\DD\setminus\{0\}$ and $f(0)=(0,x_2,0)$. Then either
$\tilde f(\DD)\subset\partial\EE$ or $\tilde f$ is a
$\tilde k_{\EE}$-extremal for
$(\tilde f(0),\tilde f(\sigma))$.
\end{lemma}

\begin{proof} We already know that $\tilde f(\DD)\subset\partial \EE$
or $\tilde f(\DD)\subset\EE$. Assume that the second property
holds. Suppose that $\tilde f$ is not a $\tilde k_{\EE}$-extremal
for $(\tilde f(0),\tilde f(\sigma))$. Then there would exist a
$g\in\OO(\DD,\EE)$ such that $g(\DD)\subset\subset\EE$,
$g(0)=\tilde f(0)$ and $g(\sigma)=\tilde f(\sigma)$. But then
$\tilde g(0)=f(0)$, $\tilde g(\sigma)=f(\sigma)$, where $\tilde
g(\lambda)=(\lambda g_1(\lambda),g_2(\lambda),\lambda
g_3(\lambda))$. But $\tilde g\in\OO(\DD,\EE)$ and $\tilde
g(\DD)\subset\subset\EE$ - a contradiction with the $\tilde
k_{\EE}$-extremality of $f$.
\end{proof}

\begin{corollary}\label{cor:3} Assume that $\phi\in\OO(\DD,\DD)$ is not an automorphism and let $\phi(0)=-C$,
where $C\in(0,1)$. Let $\omega_1,\omega_2\in\partial\DD$. Put
$f(\lambda):=(\frac{\omega_1(\phi(\lambda)+C)}{\lambda(1+C)},\omega_2\lambda\frac{1+C\phi(\lambda)}
{1+C},\omega_1\omega_2 \phi(\lambda))$. Then $f$ is a $\tilde
k_{\EE}$-extremal for $(f(0),f(\sigma))$,
$\sigma\in\DD\setminus\{0\}$.
\end{corollary}

\begin{proof} In view of Lemma~\ref{lem:1} it is sufficient to see that $f(0)=\big(\omega_1\frac{\phi^{\prime}(0)}{1+C},
0,-\omega_1\omega_2C\big)
\in\EE$. But in view of the
Schwarz-Pick Lemma (it is important here that $\phi$ is not an automorphism)
$|\phi^{\prime}(0)|<1-C^2$.

\end{proof}

\begin{rem} Note that the $\tilde k_{\EE}$-extremals from Corollary~\ref{cor:3} omit the set $S$.

Moreover, substituting $\phi\equiv-C$, $C\in[0,1)$ we get $\tilde
k_{\EE}$-extremals of the form
\begin{equation}
f(\lambda)=(0,\lambda(1-C),-C),\;\lambda\in\DD
\end{equation}
\end{rem}

\begin{proof}[Proof of Proposition~\ref{prop:1}] The proof follows directly from the fact that
$f(\lambda)=(0,\lambda(1-C),-C)$, $\lambda\in\DD$ is a $\tilde
k_{\EE}$-extremal for $(f(0),f(\lambda))$.
\end{proof}

\section{A remark on Carath\'eodory distance of the tetrablock}

Define
\begin{equation}
p_{\EE}(w,z):=\sup\{p(\Psi_{\omega}(w),\Psi_{\omega}(w)),p(\Psi_{\omega}(\sigma(w)),\Psi_{\omega}(\sigma(z))):\omega\in\partial\DD\},
\; w,z\in\EE.
\end{equation}
It is evident that $p_{\EE}\leq c_{\EE}$. We prove the following

\begin{prop}
$p_{\EE}\not\equiv c_{\EE}$.
\end{prop}
\begin{proof}
 Note that for
$\lambda\in\DD\setminus\{0\}$, $C\in(0,1)$ we get
\begin{multline}
p_{\EE}((0,0,-C),(0,\lambda(1-C),-C))=\\
\sup\{p(\omega C,\omega C+\lambda(1-C)),
p(\omega C,\frac{\omega C}{1-\omega\lambda(1-C)}):\omega\in\partial\DD\}=\\
\max\{p(0,\frac{|\lambda|}{1+C-C|\lambda|}),p(0,\frac{C|\lambda|}{1+C-|\lambda|})\}=
p(0,\frac{|\lambda|}{1+C-C|\lambda|}).
\end{multline}
Consider now the function $F(z):=\frac{z_2}{\sqrt{1+z_3-z_1z_2}}$,
$z\in\EE$. Note that $F\in\OO(\DD,\EE)$ -- to see this write
$z=\pi(A)$, where $A$ is a symmetric $2\times 2$ matrix of the
norm smaller than $1$. Then $F(z)=a_{22}/\sqrt{1-a_{12}^2}$, and
now it is sufficient to recall that $|a_{22}|^2+|a_{12}|^2<1$.

Consequently,
$c_{\EE}((0,0,-C),(0,\lambda(1-C),-C))\geq
p(0,|\lambda|\sqrt{1-C})$ but the last number is for
sufficiently small $C\in(0,1)$ larger than
$p(0,\frac{|\lambda|}{1+C-C|\lambda|})$ for small
$\lambda$.
\end{proof}

Recall that the function $p_{\mathbb G_n}$ -- similar in the
construction to $p_{\EE}$ was used while studying the problem of
the Lempert theorem for $\GG_n$. Recall that the following
(in)equalities hold: $p_{\GG_2}=c_{\GG_2}=\tilde k_{\GG_2}$ (see
\cite{Cos 2004}, \cite{Agl-You 2004}) and $p_{\GG_n}\leq
c_{\GG_n}$ and $p_{\GG_n}\not\equiv c_{\GG_n}$, $n\geq 3$ (see
\cite{NPTZ 2008}). The function $p_D$ ($D=\EE$ or $\GG_n$) or
its some generalization to more general domains may have some
connection with the class of magic functions as defined and
considered in \cite{Agl-You 2007}.

\bibliographystyle{amsplain}

\end{document}